\input amstex
\documentstyle{amsppt}
%
\catcode`@=11
\redefine\output@{%
  \def\break{\penalty-\@M}\let\par\endgraf
  \ifodd\pageno\global\hoffset=105pt\else\global\hoffset=8pt\fi  
  \shipout\vbox{%
    \ifplain@
      \let\makeheadline\relax \let\makefootline\relax
    \else
      \iffirstpage@ \global\firstpage@false
        \let\rightheadline\frheadline
        \let\leftheadline\flheadline
      \else
        \ifrunheads@ 
        \else \let\makeheadline\relax
        \fi
      \fi
    \fi
    \makeheadline \pagebody \makefootline}%
  \advancepageno \ifnum\outputpenalty>-\@MM\else\dosupereject\fi
}
\def\Beta{\mathchar"0\hexnumber@\rmfam 42}
\catcode`\@=\active
\nopagenumbers
\chardef\textvolna='176
\def\negskp{\hskip -2pt}

\chardef\degree="5E
\def\blue#1{#1}

\catcode`#=11\def\diez{#}\catcode`#=6
\catcode`&=11\catcode`&=4
\catcode`_=11\def\podcherkivanie{_}\catcode`_=8
\catcode`~=11\def\volna{~}\catcode`~=\active
\def\mycite#1{\cite{\blue{#1}}\immediate\special{ps:
     ShrHPSdict begin /ShrBORDERthickness 0 def}}
\def\myciterange#1#2#3#4{\cite{\blue{#2#3#4}}\immediate\special{ps:
     ShrHPSdict begin /ShrBORDERthickness 0 def}}
\def\mytag#1{%
    \tag#1}
\def\mythetag#1{\thetag{\blue{#1}}\immediate\special{ps:
     ShrHPSdict begin /ShrBORDERthickness 0 def}}
\def\myrefno#1{\no#1}
\def\myhref#1#2{\blue{#2}\immediate\special{ps:
     ShrHPSdict begin /ShrBORDERthickness 0 def}}
\def\myEarXivlink{\myhref{http://arXiv.org}{http:/\negskp/arXiv.org}}

\def\mytheorem#1{\csname proclaim\endcsname{Theorem #1}}
\def\mytheoremwithtitle#1#2{\csname proclaim\endcsname{Theorem #1#2}}

\def\mylemma#1{\csname proclaim\endcsname{Lemma #1}}
\def\mylemmawithtitle#1#2{\csname proclaim\endcsname{Lemma #1#2}}

\def\mycorollary#1{\csname proclaim\endcsname{Corollary #1}}

\def\myconjecture#1{\csname proclaim\endcsname{Conjecture #1}}
\def\myconjecturewithtitle#1#2{\csname proclaim\endcsname{Conjecture #1#2}}

\def\myproblem#1{\csname proclaim\endcsname{Problem #1}}
\def\myproblemwithtitle#1#2{\csname proclaim\endcsname{Problem #1#2}}
\def\mytheproblem#1{\blue{#1}\immediate\special{ps:
     ShrHPSdict begin /ShrBORDERthickness 0 def}}


\pagewidth{360pt}
\pageheight{606pt}
\topmatter
\title
On singularities of the inverse problems associated with perfect cuboids.
\endtitle
\rightheadtext{On singularities of the inverse problems \dots}
\author
John Ramsden, Ruslan Sharipov
\endauthor
\address CSR Plc, Cambridge Business Park, Cambridge, CB4 0WZ, UK
\endaddress
\email\myhref{mailto:jhnrmsdn\@yahoo.co.uk}{jhnrmsdn\@yahoo.co.uk}
\endemail
\address Bashkir State University, 32 Zaki Validi street, 450074 Ufa, Russia
\endaddress
\email\myhref{mailto:r-sharipov\@mail.ru}{r-sharipov\@mail.ru}
\endemail
\abstract
    Two cubic equations and three auxiliary equations for edges and face
diagonals of a rational perfect cuboid have been recently derived. They 
constitute a background for two inverse problems. The coefficients of cubic
equations and the right hand sides of auxiliary equations are rational
functions of two rational parameters, i\.\,e\. they have denominators. Hence
the inverse problems have singular points. These singular points are studied
in the present paper.
\endabstract
\subjclassyear{2000}
\subjclass 11D25, 11D72, 14H50, 14E05\endsubjclass
\endtopmatter
\TagsOnRight
\document

%
%
\head
1. Introduction.
\endhead
     A rational perfect cuboid is a rectangular parallelepiped whose edges and 
face diagonals are rational numbers and whose space diagonal is equal to unity:
$L=1$. Finding such a cuboid is equivalent to finding a cuboid with all integer 
edges and diagonals, which is an unsolved problem for many years (see 
\myciterange{1}{1}{--}{44}).\par
     Let $x_1$, $x_2$, $x_3$ be edges of a cuboid and let $d_1$, $d_2$, $d_3$ 
be its face diagonals. Then $x_1$, $x_2$, $x_3$ are roots of the cubic equation
$$
\hskip -2em
x^3-E_{10}\,x^2+E_{20}\,x-E_{30}=0.
\mytag{1.1}
$$
Similarly, $d_1$, $d_2$, $d_3$ are roots of the other cubic equation 
$$
\hskip -2em
d^{\kern 1pt 3}-E_{01}\,d^{\kern 1pt 2}+E_{02}\,d-E_{03}=0.
\mytag{1.2}
$$
Apart from \mythetag{1.1} and \mythetag{1.2}, the rational numbers $x_1$, $x_2$, 
$x_3$ and $d_1$, $d_2$, $d_3$ should obey the following three auxiliary equations:
$$
\hskip -2em
\aligned
&x_1\,x_2\,d_3+x_2\,x_3\,d_1+x_3\,x_1\,d_2=E_{21},\\
&x_1\,d_2+d_1\,x_2+x_2\,d_3+d_2\,x_3+x_3\,d_1+d_3\,x_1=E_{11},\\
&x_1\,d_2\,d_3+x_2\,d_3\,d_1+x_3\,d_1\,d_2=E_{12}.
\endaligned
\mytag{1.3}
$$
The cubic equations \mythetag{1.1}, \mythetag{1.2} and the auxiliary equations 
\mythetag{1.3} were obtained as a result of the series of papers \myciterange{45}
{45}{--}{50}). The coefficients $E_{10}$, $E_{20}$, $E_{30}$, $E_{01}$, $E_{02}$, 
$E_{03}$ in \mythetag{1.1} and \mythetag{1.2} as well as the right hand sides
$E_{21}$, $E_{11}$, $E_{12}$ in \mythetag{1.3} are given by explicit formulas. 
Here is the formula for $E_{11}$:
$$
\hskip -2em
E_{11}=-\frac{b\,(c^2+2-4\,c)}{b^2\,c^2+2\,b^2-3\,b^2\,c+c-b\,c^2\,+2\,b}.
\mytag{1.4}
$$
The formulas for $E_{10}$, $E_{01}$ are similar to the formula \mythetag{1.4}
for $E_{11}$:
$$
\gather
\hskip -2em
E_{10}=-\frac{b^2\,c^2+2\,b^2-3\,b^2\,c\,-c}{b^2\,c^2+2\,b^2-3\,b^2\,c
+c-b\,c^2+2\,b},
\mytag{1.5}\\
\vspace{1ex}
\hskip -2em
E_{01}=-\frac{b\,(c^2+2-2\,c)}{b^2\,c^2+2\,b^2-3\,b^2\,c+c-b\,c^2+2\,b}.
\mytag{1.6}
\endgather
$$
Below are the formulas for $E_{20}$, $E_{02}$, $E_{30}$, 
$E_{03}$, $E_{21}$, $E_{12}$ in \mythetag{1.1}, \mythetag{1.2}, 
and \mythetag{1.3}:
$$
\gather
\hskip -2em
\gathered
E_{20}=\frac{b}{2}\,(b\,c^2-2\,c-2\,b)\,(2\,b\,c^2-c^2-6\,b\,c+2
+4\,b)\,\times\\
\times\,(b\,c-1-b)^{-2}\,(b\,c-c-2\,b)^{-2},
\endgathered\qquad\quad
\mytag{1.7}\\
\vspace{1ex}
\hskip -2em
\gathered
E_{02}=\frac{1}{2}\,(28\,b^2\,c^2-16\,b^2\,c-2\,c^2-4\,b^2-b^2\,c^4
+4\,b^3\,c^4-12\,b^3\,c^3\,+\\
+\,4\,b\,c^3+24\,b^3\,c-8\,b\,c-2\,b^4\,c^4+12\,b^4\,c^3-26\,b^4\,c^2
-8\,b^2\,c^3\,+\\
+24\,b^4\,c-16\,b^3-8\,b^4)\,(b\,c-1-b)^{-2}\,(b\,c-c-2\,b)^{-2},
\endgathered\qquad\quad
\mytag{1.8}\\
\vspace{1ex}
\hskip -2em
\gathered
E_{30}=c\,b^2\,(1-c)\,(c-2)\,(b\,c^2-4\,b\,c+2+4\,b)
\,(2\,b\,c^2-c^2-4\,b\,c\,+\\
+\,2\,b)\,(b^2\,c^4-6\,b^2\,c^3+13\,b^2\,c^2-12\,b^2\,c+4\,b^2
+c^2)^{-1}\,\times\\
\times\,(b\,c-1-b)^{-2}\,(-c+b\,c-2\,b)^{-2},
\endgathered\qquad\quad
\mytag{1.9}\\
\vspace{1ex}
\hskip -2em
\gathered
E_{03}=\frac{b}{2}\,(b^2\,c^4-5\,b^2\,c^3+10\,b^2\,c^2-10\,b^2\,c+4\,b^2
+2\,b\,c+2\,c^2\,-\\
-\,b\,c^3)\,(2\,b^2\,c^4-12\,b^2\,c^3+26\,b^2\,c^2-24\,b^2\,c
+\,8\,b^2-c^4\,b+3\,b\,c^3\,-\\
-\,6\,b\,c+4\,b+c^3-2\,c^2+2\,c)\,(b^2\,c^4-6\,b^2\,c^3+13\,b^2\,c^2\,-\\
-12\,b^2\,c+4\,b^2+c^2)^{-1}\,(b\,c-1-b)^{-2}\,(-c+b\,c-2\,b)^{-2},
\endgathered\qquad\quad
\mytag{1.10}\\
\vspace{1ex}
\gathered
E_{21}=\frac{b}{2}\,(5\,c^6\,b-2\,c^6\,b^2+52\,c^5\,b^2-16\,c^5\,b
-2\,c^7\,b^2+2\,b^4\,c^8\,-\\
-\,26\,b^4\,c^7-426\,b^4\,c^5-61\,b^3\,c^6+100\,b^3\,c^5
+14\,c^7\,b^3-c^8\,b^3-20\,b\,c^2\,-\\
-\,8\,b^2\,c^2-16\,b^2\,c-128\,b^2\,c^4-200\,b^3\,c^3
+244\,b^3\,c^2+32\,b\,c^3\,+\\
+\,768\,b^4\,c^4-852\,b^4\,c^3+568\,b^4\,c^2+104\,b^2\,c^3-208\,b^4\,c
+8\,c^4\,+\\
+16\,b^3-112\,b^3\,c+142\,b^4\,c^6
+32\,b^4-2\,c^5)\,(b^2\,c^4-6\,b^2\,c^3+13\,b^2\,c^2\,-\\
-12\,b^2\,c-4\,c^3+4\,b^2+c^2)^{-1}\,(b\,c-1-b)^{-2}\,(b\,c-c-2\,b)^{-2},
\endgathered\qquad\quad
\mytag{1.11}\\
\vspace{2ex}
\hskip -2em
\gathered
E_{12}=(16\,b^6+32\,b^5-6\,c^5\,b^2+2\,c^5\,b-62\,b^5\,c^6
+62\,b^6\,c^6+16\,b^4\,-\\
-\,180\,b^6\,c^5-c^7\,b^3+18\,b^5\,c^7-12\,b^6\,c^7-2\,b^5\,c^8
+b^6\,c^8+248\,b^5\,c^2\,+\\
+\,248\,b^6\,c^2-96\,b^6\,c+321\,b^6\,c^4-180\,b^5\,c^3-144\,b^5\,c
-360\,b^6\,c^3\,+\\
+\,b^4\,c^8+8\,b^4\,c^6-6\,b^4\,c^7+18\,b^4\,c^5+7\,b^3\,c^6
+90\,b^5\,c^5-14\,b^3\,c^5\,+\\
+\,17\,b^2\,c^4+32\,b^4\,c^2+28\,b^3\,c^3-28\,b^3\,c^2-4\,b\,c^3+8\,b^3\,c
-57\,b^4\,c^4\,+\\
+\,36\,b^4\,c^3-12\,b^2\,c^3-48\,b^4\,c-c^4)\,(b^2\,c^4-6\,b^2\,c^3
+13\,b^2\,c^2\,-\\
-\,12\,b^2\,c+4\,b^2+c^2)^{-1}\,(b\,c-1-b)^{-2}\,(b\,c-c-2\,b)^{-2}.
\endgathered\qquad
\mytag{1.12}
\endgather	
$$
Based on the equations \mythetag{1.1}, \mythetag{1.2}, \mythetag{1.3}
and on the formulas \mythetag{1.4}, \mythetag{1.5}, \mythetag{1.6},
\mythetag{1.7}, \mythetag{1.8}, \mythetag{1.9}, \mythetag{1.10}, 
\mythetag{1.11}, \mythetag{1.12}, in \mycite{50} the following two 
inverse problems were formulated. 
\myproblem{1.1} Find all pairs of rational numbers $b$ and $c$ for which the
cubic equations \mythetag{1.1} and \mythetag{1.2} with the coefficients given
by the formulas \mythetag{1.5}, \mythetag{1.7},	\mythetag{1.9}, \mythetag{1.6}, 
\mythetag{1.8},	\mythetag{1.10} possess positive rational roots $x_1$, $x_2$, 
$x_3$, $d_1$, $d_2$, $d_3$ obeying the auxiliary polynomial equations 
\mythetag{1.3} whose right hand sides are given by the formulas \mythetag{1.4}, 
\mythetag{1.11}, \mythetag{1.12}. 
\endproclaim
\myproblem{1.2} Find at least one pair of rational numbers $b$ and $c$ for which 
the cubic equations \mythetag{1.1} and \mythetag{1.2} with the coefficients given
by the formulas \mythetag{1.5}, \mythetag{1.7},	\mythetag{1.9}, \mythetag{1.6}, 
\mythetag{1.8},	\mythetag{1.10} possess positive rational roots $x_1$, $x_2$, 
$x_3$, $d_1$, $d_2$, $d_3$ obeying the auxiliary polynomial equations 
\mythetag{1.3} whose right hand sides are given by the formulas \mythetag{1.4}, 
\mythetag{1.11}, \mythetag{1.12}. 
\endproclaim
     The formulas \mythetag{1.4} through \mythetag{1.12} possess denominators.
Therefore the inverse problems are singular for some values of $b$ and $c$. 
The main goal of the present paper is to study these singularities.
\head
2. The common denominator and its reduction. 
\endhead
     Let's calculate the common denominator of the fractions \mythetag{1.4} 
through \mythetag{1.12}: 
$$
\hskip -2em
\gathered
(b^2\,c^4-6\,b^2\,c^3+13\,b^2\,c^2-12\,b^2\,c+4\,b^2+c^2)\,(b\,c-1-b)^2\,\times\\
\times\,(b\,c-c-2\,b)^2\,(b^2\,c^2+2\,b^2-3\,b^2\,c+c-b\,c^2+2\,b)=0.
\endgathered
\mytag{2.1}
$$
The vanishing condition \mythetag{2.1} determines all singular points of the 
inverse problems~\mytheproblem{1.1} and \mytheproblem{1.2}. The last multiplicand
in the left hand side of the formula \mythetag{2.1} is taken from the denominators 
of \mythetag{1.4}, \mythetag{1.5}, and \mythetag{1.6}. Studying this multiplicand
we find that it is factorable. It factors as follows:
$$
b^2\,c^2+2\,b^2-3\,b^2\,c+c-b\,c^2+2\,b=(b\,c-1-b)\,(b\,c-c-2\,b).
$$
Applying this formula to \mythetag{2.1}, we reduce the equality \mythetag{2.1} to
$$
\hskip -2em
\gathered
(b^2\,c^4-6\,b^2\,c^3+13\,b^2\,c^2-12\,b^2\,c+4\,b^2+c^2)\,\times\\
\times\,(b\,c-1-b)^3\,(b\,c-c-2\,b)^3=0.
\endgathered
\mytag{2.2}
$$
Some multiplicands in \mythetag{2.1} are raised to the third power. But the exponents 
do not really affect the vanishing condition \mythetag{2.2}. Therefore we reduce 
it to
$$
\hskip -2em
\gathered
(b^2\,c^4-6\,b^2\,c^3+13\,b^2\,c^2-12\,b^2\,c+4\,b^2+c^2)\,\times\\
\times\,(b\,c-1-b)\,(b\,c-c-2\,b)=0.
\endgathered
\mytag{2.3}
$$
The left hand side of the formula \mythetag{2.3} is the reduced common denominator 
of the fractions \mythetag{1.4} through \mythetag{1.12}. It is broken into the 
product of three terms. Therefore singular points of the inverse 
problems~\mytheproblem{1.1} and \mytheproblem{1.2} are subdivided into three 
singularity subvarieties.\par
\head
3. The first singularity subvariety. 
\endhead
     The first subvariety of singular points is the most simple. It is determined
by the most simple factor in \mythetag{2.2} through the following equation: 
$$
\hskip -2em
b\,c-1-b=0.
\mytag{3.1}
$$
The equation \mythetag{3.1} is linear with respect to both $b$ and $c$. Resolving
it for $b$, we get
$$
\hskip -2em
b=\frac{1}{c-1}\text{, \ where \ }c\neq 1.
\mytag{3.2}
$$
The formula \mythetag{3.2} means that the first subvariety of singular points 
of the inverse problems~\mytheproblem{1.1} and \mytheproblem{1.2} is a rational
curve birationally equivalent to a straight line.  
\head
4. The second singularity subvariety.
\endhead
     The second subvariety of singular points is similar to the first one and is
also very simple. It is determined by the following equation:
$$
\hskip -2em
b\,c-c-2\,b=0.
\mytag{4.1}
$$
Like in \mythetag{3.2}, resolving the equation \mythetag{4.1} with respect to $b$, 
we get
$$
\hskip -2em
b=\frac{c}{c-2}\text{, \ where \ }c\neq 2.
\mytag{4.2}
$$
The formula \mythetag{4.2} means that the second subvariety of singular points 
of the inverse problems~\mytheproblem{1.1} and \mytheproblem{1.2} is a rational
curve birationally equivalent to a straight line.  
\head
5. The third singularity subvariety.
\endhead
     The third subvariety of singular points is different from the first two.
It is determined by the following equation which is quadratic with respect to 
$b$:
$$
\hskip -2em
b^2\,c^4-6\,b^2\,c^3+13\,b^2\,c^2-12\,b^2\,c+4\,b^2+c^2=0.
\mytag{5.1}
$$
The discriminant of the quadratic equation \mythetag{5.1} with respect to $b$ is
$$
\hskip -2em
D=-4\,(c-1)^2\,(c-2)^2\,c^2. 
\mytag{5.2}
$$
Looking at the discriminant formula \mythetag{5.2}, we see that we can expect 
to find rational solutions only in one of the following three cases:
$$
\xalignat 3
&\hskip -2em
c=0,
&&c=1,
&&c=2.
\mytag{5.3}
\endxalignat
$$\par
    Substituting $c=0$ into the formula \mythetag{5.1} we derive $4\,b^2=0$. 
This yields one rational point at the origin ob the $(b,c)$ plane:
$$
\pagebreak
\xalignat 2
&\hskip -2em
b=0,
&&c=0.
\mytag{5.4}
\endxalignat
$$
Substituting $c=1$ and $c=2$ into \mythetag{5.1}, we get two equalities $1=0$
and $4=0$ which are contradictory. Thus two of the three options in \mythetag{5.3}
do not actually produce any rational solutions for the equation \mythetag{5.1}.
This result is not surprising since the equation \mythetag{5.1} turns out to be 
reducible to $(c-1)^2\,(c-2)^2\,b^2+c^2=0$.\par
\head
6. Conclusions.
\endhead
     Despite being a fourth order equation in $c$, the equation \mythetag{5.1}
does not produce elliptic curves. Therefore the structure of the set of singular 
points of the inverse problems~\mytheproblem{1.1} and \mytheproblem{1.2} is
very simple. It comprises one isolated point \mythetag{5.4} and two parametric
subsets given by the formulas \mythetag{3.2} and \mythetag{4.2}. This fact can
be useful in computerized search for a solution of at least one of the 
problems~\mytheproblem{1.1} and \mytheproblem{1.2}, or in solving both of them 
if somehow it will be proved that the number of perfect cuboids is finite
and they are in a certain range.


\Refs
\ref\myrefno{1}\paper
\myhref{http://en.wikipedia.org/wiki/Euler\podcherkivanie 
brick}{Euler brick}\jour Wikipedia\publ 
Wikimedia Foundation Inc.\publaddr San Francisco, USA 
\endref
\ref\myrefno{2}\by Halcke~P.\book Deliciae mathematicae oder mathematisches 
Sinnen-Confect\publ N.~Sauer\publaddr Hamburg, Germany\yr 1719
\endref
\ref\myrefno{3}\by Saunderson~N.\book Elements of algebra, {\rm Vol. 2}\publ
Cambridge Univ\. Press\publaddr Cambridge\yr 1740 
\endref
\ref\myrefno{4}\by Euler~L.\book Vollst\"andige Anleitung zur Algebra, \rm
3 Theile\publ Kaiserliche Akademie der Wissenschaf\-ten\publaddr St\.~Petersburg
\yr 1770-1771
\endref
\ref\myrefno{5}\by Pocklington~H.~C.\paper Some Diophantine impossibilities
\jour Proc. Cambridge Phil\. Soc\. \vol 17\yr 1912\pages 108--121
\endref
\ref\myrefno{6}\by Dickson~L.~E\book History of the theory of numbers, 
{\rm Vol\. 2}: Diophantine analysis\publ Dover\publaddr New York\yr 2005
\endref
\ref\myrefno{7}\by Kraitchik~M.\paper On certain rational cuboids
\jour Scripta Math\.\vol 11\yr 1945\pages 317--326
\endref
\ref\myrefno{8}\by Kraitchik~M.\book Th\'eorie des Nombres,
{\rm Tome 3}, Analyse Diophantine et application aux cuboides 
rationelles \publ Gauthier-Villars\publaddr Paris\yr 1947
\endref
\ref\myrefno{9}\by Kraitchik~M.\paper Sur les cuboides rationelles
\jour Proc\. Int\. Congr\. Math\.\vol 2\yr 1954\publaddr Amsterdam
\pages 33--34
\endref
\ref\myrefno{10}\by Bromhead~T.~B.\paper On square sums of squares
\jour Math\. Gazette\vol 44\issue 349\yr 1960\pages 219--220
\endref
\ref\myrefno{11}\by Lal~M., Blundon~W.~J.\paper Solutions of the 
Diophantine equations $x^2+y^2=l^2$, $y^2+z^2=m^2$, $z^2+x^2
=n^2$\jour Math\. Comp\.\vol 20\yr 1966\pages 144--147
\endref
\ref\myrefno{12}\by Spohn~W.~G.\paper On the integral cuboid\jour Amer\. 
Math\. Monthly\vol 79\issue 1\pages 57-59\yr 1972 
\endref
\ref\myrefno{13}\by Spohn~W.~G.\paper On the derived cuboid\jour Canad\. 
Math\. Bull\.\vol 17\issue 4\pages 575-577\yr 1974
\endref
\ref\myrefno{14}\by Chein~E.~Z.\paper On the derived cuboid of an 
Eulerian triple\jour Canad\. Math\. Bull\.\vol 20\issue 4\yr 1977
\pages 509--510
\endref
\ref\myrefno{15}\by Leech~J.\paper The rational cuboid revisited
\jour Amer\. Math\. Monthly\vol 84\issue 7\pages 518--533\yr 1977
\moreref see also Erratum\jour Amer\. Math\. Monthly\vol 85\page 472
\yr 1978
\endref
\ref\myrefno{16}\by Leech~J.\paper Five tables relating to rational cuboids
\jour Math\. Comp\.\vol 32\yr 1978\pages 657--659
\endref
\ref\myrefno{17}\by Spohn~W.~G.\paper Table of integral cuboids and their 
generators\jour Math\. Comp\.\vol 33\yr 1979\pages 428--429
\endref
\ref\myrefno{18}\by Lagrange~J.\paper Sur le d\'eriv\'e du cuboide 
Eul\'erien\jour Canad\. Math\. Bull\.\vol 22\issue 2\yr 1979\pages 239--241
\endref
\ref\myrefno{19}\by Leech~J.\paper A remark on rational cuboids\jour Canad\. 
Math\. Bull\.\vol 24\issue 3\yr 1981\pages 377--378
\endref
\ref\myrefno{20}\by Korec~I.\paper Nonexistence of small perfect 
rational cuboid\jour Acta Math\. Univ\. Comen\.\vol 42/43\yr 1983
\pages 73--86
\endref
\ref\myrefno{21}\by Korec~I.\paper Nonexistence of small perfect 
rational cuboid II\jour Acta Math\. Univ\. Comen\.\vol 44/45\yr 1984
\pages 39--48
\endref
\ref\myrefno{22}\by Wells~D.~G.\book The Penguin dictionary of curious and 
interesting numbers\publ Penguin publishers\publaddr London\yr 1986
\endref
\ref\myrefno{23}\by Bremner~A., Guy~R.~K.\paper A dozen difficult Diophantine 
dilemmas\jour Amer\. Math\. Monthly\vol 95\issue 1\yr 1988\pages 31--36
\endref
\ref\myrefno{24}\by Bremner~A.\paper The rational cuboid and a quartic surface
\jour Rocky Mountain J\. Math\. \vol 18\issue 1\yr 1988\pages 105--121
\endref
\ref\myrefno{25}\by Colman~W.~J.~A.\paper On certain semiperfect cuboids\jour
Fibonacci Quart.\vol 26\issue 1\yr 1988\pages 54--57\moreref see also\nofrills 
\paper Some observations on the classical cuboid and its parametric solutions
\jour Fibonacci Quart\.\vol 26\issue 4\yr 1988\pages 338--343
\endref
\ref\myrefno{26}\by Korec~I.\paper Lower bounds for perfect rational cuboids 
\jour Math\. Slovaca\vol 42\issue 5\yr 1992\pages 565--582
\endref
\ref\myrefno{27}\by Guy~R.~K.\paper Is there a perfect cuboid? Four squares 
whose sums in pairs are square. Four squares whose differences are square 
\inbook Unsolved Problems in Number Theory, 2nd ed.\pages 173--181\yr 1994
\publ Springer-Verlag\publaddr New York 
\endref
\ref\myrefno{28}\by Rathbun~R.~L., Granlund~T.\paper The integer cuboid table 
with body, edge, and face type of solutions\jour Math\. Comp\.\vol 62\yr 1994
\pages 441--442
\endref
\ref\myrefno{29}\by Van Luijk~R.\book On perfect cuboids, \rm Doctoraalscriptie
\publ Mathematisch Instituut, Universiteit Utrecht\publaddr Utrecht\yr 2000
\endref
\ref\myrefno{30}\by Rathbun~R.~L., Granlund~T.\paper The classical rational 
cuboid table of Maurice Kraitchik\jour Math\. Comp\.\vol 62\yr 1994
\pages 442--443
\endref
\ref\myrefno{31}\by Peterson~B.~E., Jordan~J.~H.\paper Integer hexahedra equivalent 
to perfect boxes\jour Amer\. Math\. Monthly\vol 102\issue 1\yr 1995\pages 41--45
\endref
\ref\myrefno{32}\by Rathbun~R.~L.\paper The rational cuboid table of Maurice 
Kraitchik\jour e-print \myhref{http://arxiv.org/abs/math/0111229}{math.HO/0111229} 
in Electronic Archive \myEarXivlink
\endref
\ref\myrefno{33}\by Hartshorne~R., Van Luijk~R.\paper Non-Euclidean Pythagorean 
triples, a problem of Euler, and rational points on K3 surfaces\publ e-print 
\myhref{http://arxiv.org/abs/math/0606700}{math.NT/0606700} 
in Electronic Archive \myEarXivlink
\endref
\ref\myrefno{34}\by Waldschmidt~M.\paper Open diophantine problems\publ e-print 
\myhref{http://arxiv.org/abs/math/0312440}{math.NT/0312440} 
in Electronic Archive \myEarXivlink
\endref
\ref\myrefno{35}\by Ionascu~E.~J., Luca~F., Stanica~P.\paper Heron triangles 
with two fixed sides\publ e-print \myhref{http://arxiv.org/abs/math/0608185}
{math.NT/0608} \myhref{http://arxiv.org/abs/math/0608185}{185} in Electronic 
Archive \myEarXivlink
\endref
\ref\myrefno{36}\by Ortan~A., Quenneville-Belair~V.\paper Euler's brick
\jour Delta Epsilon, McGill Undergraduate Mathematics Journal\yr 2006\vol 1
\pages 30-33
\endref
\ref\myrefno{37}\by Knill~O.\paper Hunting for Perfect Euler Bricks\jour Harvard
College Math\. Review\yr 2008\vol 2\issue 2\page 102\moreref
see also \myhref{http://www.math.harvard.edu/\volna knill/various/eulercuboid/index.html}
{http:/\negskp/www.math.harvard.edu/\textvolna knill/various/eulercuboid/index.html}
\endref
\ref\myrefno{38}\by Sloan~N.~J.~A\paper Sequences 
\myhref{http://oeis.org/A031173}{A031173}, 
\myhref{http://oeis.org/A031174}{A031174}, and \myhref{http://oeis.org/A031175}
{A031175}\jour On-line encyclopedia of integer sequences\publ OEIS Foundation 
Inc.\publaddr Portland, USA
\endref
\ref\myrefno{39}\by Stoll~M., Testa~D.\paper The surface parametrizing cuboids
\jour e-print \myhref{http://arxiv.org/abs/1009.0388}{arXiv:1009.0388} 
in Electronic Archive \myEarXivlink
\endref
\ref\myrefno{40}\by Sharipov~R.~A.\paper A note on a perfect Euler cuboid.
\jour e-print \myhref{http://arxiv.org/abs/1104.1716}{arXiv:1104.1716} 
in Electronic Archive \myEarXivlink
\endref
\ref\myrefno{41}\by Sharipov~R.~A.\paper Perfect cuboids and irreducible 
polynomials\jour Ufa Mathematical Journal\vol 4, \issue 1\yr 2012\pages 153--160
\moreref see also e-print \myhref{http://arxiv.org/abs/1108.5348}{arXiv:1108.5348} 
in Electronic Archive \myEarXivlink
\endref
\ref\myrefno{42}\by Sharipov~R.~A.\paper A note on the first cuboid conjecture
\jour e-print \myhref{http://arxiv.org/abs/1109.2534}{arXiv:1109.2534} 
in Electronic Archive \myEarXivlink
\endref
\ref\myrefno{43}\by Sharipov~R.~A.\paper A note on the second cuboid conjecture.
Part~\uppercase\expandafter{\romannumeral 1} 
\jour e-print \myhref{http://arxiv.org/abs/1201.1229}{arXiv:1201.1229} 
in Electronic Archive \myEarXivlink
\endref
\ref\myrefno{44}\by Sharipov~R.~A.\paper A note on the third cuboid conjecture.
Part~\uppercase\expandafter{\romannumeral 1} 
\jour e-print \myhref{http://arxiv.org/abs/1203.2567}{arXiv:1203.2567} 
in Electronic Archive \myEarXivlink
\endref
\ref\myrefno{45}\by Sharipov~R.~A.\paper Perfect cuboids and multisymmetric 
polynomials\jour e-print \myhref{http://arxiv.org/abs/1203.2567}
{arXiv:1205.3135} in Electronic Archive \myEarXivlink
\endref
\ref\myrefno{46}\by Sharipov~R.~A.\paper On an ideal of multisymmetric polynomials 
associated with perfect cuboids\jour e-print \myhref{http://arxiv.org/abs/1206.6769}
{arXiv:1206.6769} in Electronic Archive \myEarXivlink
\endref
\ref\myrefno{47}\by Sharipov~R.~A.\paper On the equivalence of cuboid equations and 
their factor equations\jour e-print \myhref{http://arxiv.org/abs/1207.2102}
{arXiv:1207.2102} in Electronic Archive \myEarXivlink
\endref
\ref\myrefno{48}\by Sharipov~R.~A.\paper A biquadratic Diophantine equation associated 
with perfect cuboids\jour e-print \myhref{http://arxiv.org/abs/1207.4081}
{arXiv:1207.4081} in Electronic Archive \myEarXivlink
\endref
\ref\myrefno{49}\by Ramsden~J.~R.\paper A general rational solution of an equation 
associated with perfect cuboids\jour e-print \myhref{http://arxiv.org/abs/1207.5339}
{arXiv:1207.5339} in Electronic Archive \myEarXivlink
\endref
\ref\myrefno{50}\by Ramsden~J.~R., Sharipov~R.~A.\paper Inverse problems associated 
with perfect cuboids\jour e-print \myhref{http://arxiv.org/abs/1207.6764}
{arXiv:1207.6764} in Electronic Archive \myEarXivlink
\endref
\endRefs
\enddocument
\end